\documentclass[12pt]{article}
\usepackage{amsmath,amsfonts,amsthm}

\DeclareMathAlphabet{\mathbbold}{U}{bbold}{m}{n}

\DeclareSymbolFont{rsfscript}{OMS}{rsfs}{m}{n}
\DeclareSymbolFontAlphabet{\mathrsfs}{rsfscript}

\DeclareFontFamily{OMS}{rsfs}{\skewchar\font'177}
\DeclareFontShape{OMS}{rsfs}{m}{n}{%
      <5> rsfs5
      <6> <7> rsfs7
      <8> <9> <10> rsfs10
      <10.95> <12> <14.4> <17.28> <20.74> <24.88> rsfs10
      }{}

\DeclareMathOperator{\Hom}{Hom}
\DeclareMathOperator{\Res}{Res}
\DeclareMathOperator{\Ind}{Ind}

\DeclareMathOperator{\PF}{PF}
\DeclareMathOperator{\Fix}{Fix}
\makeatletter

\newcommand{\theoremName}{Theorem}
\newcommand{\proofName}{Proof}
\renewcommand{\proofname}{\proofName}

\theoremstyle{plain}

\newtheorem {remark}{Remark}
\newtheorem {example}{Example}
\newtheorem {proposition}{Proposition}

\let\@newpf\proof \let\proof\relax 
\newenvironment{proof}{\@newpf[\proofname]}{\qed\endtrivlist}

\begin{document}
\title{A remark on Frobenius characters\\ for set representations of symmetric groups}
\author{Vladimir Dotsenko}
\date{}

\maketitle

\begin{abstract}
For any set representation (permutation representation) of the symmetric group $S_n$, we give combinatorial interpretation for coefficients of its Frobenius character expanded in the basis of monomial symmetric functions.
\end{abstract}

\section{Introduction.}
Set representations (permutation representations) of groups are representations arising from actions on finite sets. In other words, for any group homomorphism $\phi\colon G\to S_k$ we can consider a representation of $G$ which is obtained from restriction of the tautological $k$-dimensional representation of $S_k$ to~$G$. Such representations are called set representations. From this definition, it is easy to see that, unlike the general case of group representations, the values of character for such a representation have a nice combinatorial interpretation in terms of fixed points. In the case when $G$ itself is a symmetric group, the values of the character can be arranged in a generating function called Frobenius character which is an element of the ring of symmetric polynomials. The ring of symmetric polynomials has several important bases, and it is natural to ask whether or not the expansions of the Frobenius character for a permutation representation relative to these bases have simple combinatorial interpretations as well. In this paper we discuss this question; in particular, we discuss such an interpretation for the basis of monomial symmetric polynomials; this result has much in common with the main theorem of Polya enumeration theory and also is immediate from Frobenius reciprocity point of view. This observation is quite simple, and probably belongs to some sort of mathematical folklore, though the author is not aware of suitable references. For the author, the starting motivation was his search for combinatorical interpretations of corresponding coefficients for the case of action on parking functions; it was primarily motivated by an interpretation of parking functions module in terms of lattice vertex operator algebras that he obtained recently. These results will appear elsewhere~\cite{VD}. 
 
The paper is organised as follows. In Section \ref{notation} we list some classical definitions related to combinatorics and representation theory of symmetric groups. Section \ref{frobenius} contains the combinatorial formula discussed above. Finally, in Section \ref{parking} we discuss an interesting particular case of our formula that corresponds to the action of symmetric groups on parking functions, obtaining thus a new proof for one well known formula. 

The author is grateful to Andrei Zelevinsky for the most useful comments on the overly ambitious previous version of this paper.

\section{Definitions and notation.}\label{notation}

Throughout the paper, $[n]$ denotes the set $\{1,2,\ldots,n\}$. Greek letters $\lambda$, $\mu$ etc. denote partitions, i.e. decompositions of an integer number $N$ into a sum of nonnegative numbers $N=\lambda_1+\lambda_2+\ldots+\lambda_N$; in this notation $\lambda$ denotes the sequence $\lambda_1$, \ldots, $\lambda_N$, where we assume that numbers are rearranged so that $\lambda_1\ge\lambda_2\ge\ldots\ge\lambda_N\ge0$.
A brief notation for ``$\lambda$ is a partition of $N$'' is $\lambda\vdash N$. By the definition, $m_s(\lambda)$ denotes the number of parts in $\lambda$ that are equal to~$s$. The number of nonzero parts in $\lambda$ is called length of $\lambda$ and is denoted by $\ell(\lambda)$.

Suppose that a group $G$ acts on a finite set $M$. We use the standard notation $M/G$ for the set of orbits of action of this action. The space where the corresponding set representation of $G$ is realised is the set of functions on~$M$; it is denoted by $C(M)$. The following proposition is well known.

\begin{proposition}\label{set-char}
For the character of $C(M)$, we have for any $g\in G$
 $$
\chi_{C(M)}(g)=\#\Fix(g),
 $$ 
where $\Fix(g)=M^g$ is the set of fixed points of~$g$.
\end{proposition}

To any representation $V$ of the symmetric group $S_n$, one can assign the following 
polynomial $\mathcal{F}_V$ in variables $p_1$, \ldots, $p_n$. Recall that conjugacy classes in $S_n$ are in one-to-one
correspondence with partitions of $n$: a partition $\lambda$ corresponds to the conjugacy class of permutations whose decomposition into a product of disjoint cycles contains cycles of lengths $\lambda_1$, $\lambda_2$, \ldots 
The polynomial $\mathcal{F}_V$, usually called the Frobenius character of $V$ \cite{Mac}, is given by the formula
 $$
\mathcal{F}_V(p_1,p_2,\ldots,p_n)=\sum_{\lambda\vdash n}\chi_V(\lambda)\frac{p_\lambda}{z\lambda},
 $$
where $\chi_V(\lambda)$ is the value of character on an arbitrary representative of the conjugacy class $\lambda$,  $p_\lambda=p_{\lambda_1}\cdot\ldots\cdot p_{\lambda_n}$, and $z_\lambda$ is the number of elements in the centraliser of any representative of the conjugacy class; $z_\lambda=\prod_s s^{m_s(\lambda)}m_s(\lambda)!$. If we interpret the variables $p_k$ as Newton power sums $x_1^k+\ldots+x_n^k$ in the ring of symmetric polynomials, the Frobenius character becomes an element of this ring. Under this correspondence, irreducible representations correspond to Schur symmetric functions $s_\lambda$, which form a basis for this ring over the ground field. Other frequently used bases (besides the power sums basis $p_\lambda$) are $e_\lambda=e_{\lambda_1}\cdot\ldots\cdot e_{\lambda_n}$, where $e_p$  is the elementary symmetric polynomial (the symmetrisation of $x_1x_2\cdot\ldots\cdot x_p$), $h_\lambda=h_{\lambda_1}\cdot\ldots\cdot h_{\lambda_n}$, where $h_p$ is the complete symmetric function (sum of all monomials of degree $p$), and $m_\lambda$, the monomial symmetric polynomials ($m_\lambda$ is the symmetrisation of $x_1^{\lambda_1}\cdot\ldots\cdot x_n^{\lambda_n}$).

\section{A formula for Frobenius characters.}\label{frobenius}

Frobenius characters are just generating functions for characters, so for a general representation of $S_n$ they might be quite complicated. In the case of set representations, they carry some information about the group action, as one can readily see from Proposition \ref{set-char} that expands Frobenius characters in the power sums basis in terms of what we know about the group action. The following example shows that not for all bases coefficients of expansion are nonnegative.

\begin{example}
For the only nontrivial homomorphism $\phi\colon S_4\to S_3$ (quotient modulo the Klein group), the character of the corresponding set representation is equal to 
 $$
s_4+s_{2,2}=h_4-h_{3,1}+h_{2,2}=-e_4+e_{3,1}+2e_{2,2}-3e_{2,1,1}+e_{1,1,1,1}.
 $$ 
Thus, for the expansion relative to the $h$-basis or the $e$-basis the positivity property for coefficients does not hold, so we should not expect immediate combinatorial interpretations. 
\end{example}

For the case of the Schur polynomials, the corresponding coefficients are positive and are just multiplicities of
irreducibles, and no immediate combinatorial description seems to be known. Another case when the coefficients are positive (for obvious reasons) and so we might expect some combinatorics behind is the case of monomial symmetric polynomials. The following proposition establishes a simple interpretation of coefficients in term of the group action.

\begin{proposition}
Assume that $V=C(M)$ is a set-representation of $S_n$. Then we have the following expansion for the Frobenius character of $V$ in the basis of monomial symmetric functions:
 $$
\mathcal{F}_{V}=\sum_{\mu\vdash n}\#(M/S_\mu) m_\mu,
 $$
where $S_\mu=S_{\mu_1}\times S_{\mu_2}\times\dots\times S_{\mu_k}$.
\end{proposition}

\begin{remark}
This formula should remind reader of the key formula of the Polya enumeration theory~\cite{Stan}. The proof below, exactly like the proof of the main result of Polya theory, makes use of only Burnside's formula, and hopefully will lead the reader to a better understanding of Polya theory as well. There is also a following more simple representation-theoretical explanation of the above fact (communicated to me by A.Zelevinsky): the coefficient of $m_\mu$ in $\mathcal{F}_V$ is (for the standard inner product $\langle\cdot,\cdot\rangle$ on the ring of symmetric polynomials
 $$
\langle \mathcal{F}_V, h_\mu\rangle=\dim\Hom_{S_n}(V,\Ind_{S_\mu}^{S_n}\mathbbold{1})=\dim\Hom_{S_\mu}(\Res_{S_n}^{S_\mu}V,\mathbbold{1}),  
 $$
where the last equality uses Frobenius reciprocity theorem. The multiplicuty of trivial representation of $S_\mu$ in $\Res_{S_n}^{S_\mu}V$
is obviously equal to the number of orbits of $S_\mu$ in $M$.
\end{remark}

\begin{proof}
From Proposition \ref{set-char} we have the following formula for the Frobenius character in the power sums basis:
 $$
\mathcal{F}_{V}=\sum_{\lambda\vdash n}\#\Fix(\lambda) \frac{p_\lambda}{z_\lambda}.
 $$
To rewrite this in the monomial basis, we use the following
\begin{proposition}[{\cite[Ex.~7.14]{Mac}}]
In the ring of symmetric polynomials
 $$
p_\lambda=\sum_{\mu\vdash n}L_{\lambda\mu}m_\mu,
 $$
where $L_{\lambda\mu}$ is equal to the number of mappings $\phi$ from $[\ell(\lambda)]$ to $[\ell(\mu)]$ for which $\sum_{\phi(i)=j}\lambda_i=\mu_j$ for all $j=1,\ldots,\ell(\mu)$. 
\end{proposition}
\begin{proof}
Indeed, 
 $$
p_\lambda=\sum_i x_i^{\lambda_1}\sum_i x_i^{\lambda_2}\cdot\ldots\cdot\sum_i x_i^{\lambda_{\ell(\lambda)}},
 $$
so monomials that occur in the expansion of $p_\lambda$ are 
 $$
x_{i_1}^{\lambda_1}x_{i_2}^{\lambda_2}\cdot\ldots\cdot x_{i_{\ell(\lambda)}}^{\lambda_{\ell(\lambda)}},
 $$ 
where some of subscripts $i_1$, $i_2$,\ldots, $i_{\ell(\lambda)}$ might be equal to each other, which is exactly what we want.
\end{proof}
This leads to the following equivalent form of our formula:
 $$
\mathcal{F}_{V}=\sum_{\lambda\vdash n}\frac{\#\Fix(\lambda)}{z_\lambda}\sum_{\mu\vdash n}L_{\lambda\mu}m_\mu=
\sum_{\mu\vdash n}m_\mu\sum_{\lambda\vdash n}\frac{\#\Fix(\lambda)L_{\lambda\mu}}{z_\lambda}.
 $$
This means that we need to prove that for all $\mu$
 $$
\sum_{\lambda\vdash n}\frac{\#\Fix(\lambda)L_{\lambda\mu}}{z_\lambda}=\# M/S_\mu.
 $$
To compute the number of orbits, we use the following well known fact.
\begin{proposition}[Burnside's formula~\cite{Stan}]
If a group $G$ acts on a finite set $M$, then the number of orbits is equal to the average number of fixed points:
 $$
\# M/G=\frac{1}{\#G}\sum_{g\in G}\#\Fix(g).
 $$
\end{proposition}
Applying of Burnside's formula to the case of the group $S_\mu$ acting on $M$, we can rewrite our formula as
\begin{equation}\label{form1}
\sum_{\lambda\vdash n} \frac{\#\Fix(\lambda)L_{\lambda\mu}}{z_\lambda}=\\=\frac1{\mu_1!\mu_2!\ldots\mu_{\ell(\mu)}!}
\sum_{g\in S_\mu}\#\Fix(g). 
\end{equation}
Now the proof can be finished in a quite straightforward way. First, we group the summands in the right-hand side that correspond to elements $g$ from the same conjugacy class. The number of occurences of such a summand is equal to the cardinality of the corresponding conjugacy class. Conjugacy classes in $S_\mu$ are in one-to-one correspondence with sequences of $\ell(\mu)$ partitions with $|\mu^{(i)}|=\mu_i$ for all~$i=1,\ldots,\ell(\mu)$. Denote such a sequence by $(\mu^{(1)},\mu^{(2)},\ldots,\mu^{({\ell(\mu)})})$. The cardinality of the corresponding conjugacy class is equal to the index of the centraliser of an element from the conjugacy class, which is equal to 
 $$
\frac{\mu_1!\mu_2!\ldots\mu_l!}{z_{\mu^{(1)}}z_{\mu^{(2)}}\ldots z_{\mu^{({\ell(\mu)})}}}.
 $$
Thus, we managed to rewrite our formula as
\begin{multline}\label{form2}
\sum_{\lambda\vdash n}\frac{\#\Fix(\lambda)L_{\lambda\mu}}{z_\lambda}=\\=
\sum_{|\mu^{(i)}|=\mu_i}
\frac{1}{z_{\mu^{(1)}}z_{\mu^{(2)}}\ldots z_{\mu^{({\ell(\mu)})}}}\#\Fix((\mu^{(1)},\mu^{(2)},\ldots,\mu^{({\ell(\mu)})})). 
\end{multline}
Some of conjugacy classes in $S_\mu$ correspond to the same conjugacy class in $S_n$. More precisely, occurrences of the conjugacy class $\lambda$ in the right hand side of \eqref{form2} are numbered by different ways to distribute parts of $\lambda$ between ${\ell(\mu)}$ partitions $\eta^{(1)}$, $\eta^{(2)}$, \ldots, $\eta^{({\ell(\mu)})}$. For such a way, the contribution to the coefficient of $\#\Fix(\lambda)$ is 
 $$
\frac{1}{z_{\eta^{(1)}}z_{\eta^{(2)}}\ldots z_{\eta^{({\ell(\mu)})}}}=\frac{1}{z_\lambda}
\prod_s\binom{m_s(\lambda)}{m_s(\eta^{(1)}),\ldots,m_s(\eta^{({\ell(\mu)})})}.
 $$
To complete the proof, it remains to notice that the following formula is obvious from the definition of numbers $L_{\lambda\mu}$:
 $$
L_{\lambda\mu}=\sum_{\mu^{(i)}\vdash\mu_i}\prod_s\binom{m_s(\lambda)}{m_s(\mu^{(1)}),\ldots,m_s(\mu^{({\ell(\mu)})})}.
 $$
\end{proof}

\section{An example: parking functions module.}\label{parking}

A parking function of length $n$ is a function $f\colon[n]\to[n]$ satisfying the condition
$\#f^{-1}([k])\ge k$ for all $1\le k\le n$. The name ``parking function'' has the following
combinatorial explanation. Imagine a one-way street with parking spaces labeled from $1$ to $n$.
There are $n$ cars which want to park along the street, and each car $i$ has a preferred parking
space $f(i)$. The cars arrive successively at the head of the street; a car drives directly to its preferred
parking space. If the space is not occupied, the car parks there; otherwise it continues to the next
unoccupied space. If any car reaches the end of the street without having parked, the process fails.
In this terms, a parking function is a parking preference for this process to succeed.

The set of all parking functions of length~$n$ is denoted
by $\PF_n$. The parking function condition is stable under permutations of arguments, so $S_n$ acts on
$\PF_n$. The corresponding permutation representation have recently got a very interesting and unexpected
interpretation via ``diagonal harmonics''~\cite{Hai}. 

Here we use our theorem to derive a combinatorial proof for the following proposition.

\begin{proposition}[{\cite[Prop.~2.2]{Stan1}}] We have
 $$
\mathcal{F}_{\PF_n}=\sum_{\mu\vdash n}\frac{1}{n+1}\left[\prod_{i}\binom{\mu_i+n}{n}\right]m_\mu.
 $$
\end{proposition}
\begin{proof}
We use our main theorem to see that it is enough to prove that coefficients in the right hand side count orbits of the corresponding subgroups. To prove that, we will use a remarkable proof of the formula $\#\PF_n=(n+1)^{n-1}$ due to Pollak, \cite[p.~13]{FR}. Note that this formula is a particular case of our statement for $\lambda=(1,1,\ldots,1)$. Indeed, that partition corresponds to the trivial subgroup of $S_n$ (so the number of orbits is equal to the number of parking functions). 

The Pollak's proof goes as follows. If we consider the parking process on a circular one-way street with $n+1$ parking spaces $1$, \ldots, $n+1$, then the process can't fail, and exactly one parking space will remain unoccupied. Among the $n+1$ distinct rotations of any given preference function $f\colon[n]\to[n+1]$ precisely one is a parking function --- the one for which the unoccupied space is $n+1$.

Let us now prove the general statement. An $S_\mu$-orbit in $\PF_n$ is completely determined by the preference set of the first $\mu_1$ cars, the preference set of the next $\mu_2$ cars etc. Again, considering the process for a circular road, first want to compute the total number of preferences. It is equal to the number of ways to define multisets of $\mu_1$, $\mu_2$, \ldots elements, using numbers from $1$ to $n+1$ as elements. For a multiset of $\mu_i$ elements, the number of ways is $\binom{\mu_i+n}{\mu_i}=\binom{\mu_i+n}{n}$, and the total number of preferences is equal to the product of these numbers. Among the rotations of any given preference only one gives a parking function, so to compute the number of orbits we should divide this product by $(n+1)$, which completes the proof.
\end{proof}

{\footnotesize

\medskip

\noindent 
School of Mathematcs, Trinity College Dublin, Dublin 2, Ireland\\
\texttt{vdots@maths.tcd.ie}}

\end{document}